\documentclass[12pt]{amsart}

\def\nz{\ifmmode {I\hskip -3pt N} \else {\hbox {$I\hskip -3pt N$}}\fi}
\def\zz{\ifmmode {Z\hskip -4.8pt Z} \else
       {\hbox {$Z\hskip -4.8pt Z$}}\fi}
\def\qz{\ifmmode {Q\hskip -5.0pt\vrule height6.0pt depth 0pt
       \hskip 6pt} \'else {\hbox
       {$Q\hskip -5.0pt\vrule height6.0pt depth 0pt\hskip 6pt$}}\fi}
\def\rz{\ifmmode {I\hskip -3pt R} \else {\hbox {$I\hskip -3pt R$}}\fi}
\def\cz{\ifmmode {C\hskip -4.8pt\vrule height5.8pt\hskip 6.3pt} \else
       {\hbox {$C\hskip -4.8pt\vrule height5.8pt\hskip 6.3pt$}}\fi}

\def\qed{\hbox {\hskip 1pt \vrule width 4pt height 6pt depth 1.5pt
        \hskip 1pt}\\}% cqfd
\def\tr {{\rm \;tr \;}}

\newcommand {\pa}{\partial}
\newcommand {\ar}{\rightarrow}

\newcommand{\vD}{\mathcal D}
\newcommand{\bR}{\mathbb  R}

\newcommand{\bN}{\mathbb  N}

\makeatletter
\@addtoreset{equation}{section}

\makeatother

\newtheorem{theorem}{Theorem}[section]
\newtheorem{lemma}[theorem]{Lemma}
\newtheorem{proposition}[theorem]{Proposition}

\newtheorem{remark}[theorem]{Remark}
\newtheorem{corollary}[theorem]{Corollary}

\begin{document}

\title[Semiclassical analysis of a nonlinear eigenvalue problem]{Semiclassical analysis of a nonlinear eigenvalue problem and
non analytic hypoellipticity}

\author{Bernard HELFFER}
\address{\hskip-\parindent
Bernard HELFFER\\
D\'epartement de Math\'ematiques, UMR CNRS 8628\\
Universit\'e Paris-Sud, Bat. 425 \\
91405 Orsay Cedex,  France        }
%\email{email 1}

\author{Didier ROBERT}
\address{\hskip-\parindent
Didier ROBERT \\
Laboratoire de Math\'ematiques Jean Leray\\
UMR CNRS 6629\\
D\'epar\-tement de Math\'e\-ma\-tiques\\ 
Universit\'e de Nantes \\
44322 Nantes Cedex 3, France}
%\email{email 2}

\author{Xue Ping WANG}
\address{\hskip-\parindent
 Xue Ping WANG \\
Laboratoire de Mat\-h\'e\-ma\-tiques Jean Leray\\
UMR CNRS 6629\\
D\'epartement de Math\'ematiques\\ 
Universit\'e de Nantes \\
44322 Nantes Cedex 3, France}

\thanks{This work was achieved when two of the authors (B. H. and X.P. W.) were at the MSRI in Berkeley.
They would like to thank the Institute and the organizers of special semester {\it ``Semi-classical Analysis''} for invitation. The  authors would like to thank S.~Chanillo, M.~Christ, A.~Laptev and F.~Tr\`eves for useful discussions, comments  or correspondence, and J.~Sj\"ostrand for informing us about the reference \cite{Ne}. 
B. H. has  a partial support by the SPECT ESF european programme. 
}

\maketitle

\begin{abstract}
We give a semiclassical analysis of a nonlinear eigenvalue problem
arising from the study of the failure of analytic hypoellipticity and
obtain  a general  family of  hypoelliptic,
but  not analytic hypoelliptic operators.
\end{abstract}

\section{Introduction}
We are interested in a family of operators  of the type
\begin{equation}
H(x, D_x, \lambda) = - \Delta + (P(x) - \lambda )^2\;,
\end{equation}
where $x\mapsto P(x)$ is a polynomial of degree $m$.
In the study of the failure of analytic hypoellipticity, one approach 
consists in showing  that
the nonlinear eigenvalue problem
\begin{equation} \label{p1}
H(x, D_x, \lambda)v=0
\end{equation}
has at least one non trivial solution $(\lambda, v)\in C \times \mathcal
S(\mathbb R^n)$
($v\neq 0$).
This has been used by many authors including Baouendi-Goulaouic, 
Helffer, Christ \cite{Chr1, Chr2, Chr3, Chr4}, Hanges-Himonas \cite{HH},  Chanillo \cite{Cha}  and quite 
recently by 
Chanillo-Helffer-Laptev \cite{CHL}, where the reader can find a more 
extensive list of references.
All the results lead  to the formulation of
 a conjecture
 by Tr\`eves \cite{Tr} giving a necessary and sufficient condition of 
analytic hypoellipticity extending \cite{Met}.
It could be a rather natural conjecture that,
 when  $x\mapsto P(x)$ is an homogeneous elliptic polynomial on 
$\mathbb  R^n$ of order $m>1$, (\ref{p1}) has at least one non-trivial 
solution.
 This result is proved  (\cite{Chr1,PhRo})  for $n=1$  and for $m \ge 2n$
 when $n=2,3$ (\cite{CHL}).  Our aim is to  provide a semiclassical  approach to this 
problem.  Our analysis concerns actually a more general class of operators of the form
\begin{equation}\label{PQ}
 \sum_{j=1}^p  D_{x_j}^2  +
 \left(P (x_1,\ldots,x_p)  D_{x_{p+1}} - D_{x_{p+2}}\right)^2 + \left(Q (x_1,\ldots,x_p)  D_{x_{p+1}}\right)^2\;.
\end{equation}
When specialized to the case (\ref{p1}), we obtain the following

 \begin{theorem}\label{neven}.\\
 Let $n$ be even.  Let $P$ be any real polynomial of  degree $ m \ge 2$
such that the homogeneous part of degree $m$ is elliptic.
 Then the  nonlinear eigenvalue problem  (\ref{p1}) has at least one non  trivial solution.
 \end{theorem}

Its now standard corollary to analytic hypoellipticity (see for example 
\cite{CHL}) is~:

\begin{corollary}.\\
Let $P$ be elliptic  homogeneous of degree $m\geq 2$ on $\mathbb R^p$, 
with an even $p >0$, then
the operator  on $\mathbb R^{p+2}$~:
$$
P(x,D_x):= \sum_{j=1}^p  D_{x_j}^2  +
 \left(P (x_1,\ldots,x_p)  D_{x_{p+1}} - D_{x_{p+2}}\right)^2\;,
$$
is not  analytic hypoelliptic at $0$.
\end{corollary}

\begin{remark}.~\\
As we shall discuss more deeply in Section \ref{Section5}, G.~M\'etivier
 already showed in \cite{Met4}, that $P(x,D_x)$ is not hypoanalytic
 in any open set. This is indeed a sum of squares with an odd number
 of linearly independent vector fields. But the corollary presented here
 and its proof give  a  stronger information at $0$.
\end{remark}

 Other examples in any dimension related to (\ref{PQ}) are also given in Section \ref{Section4}.
 These examples give an additional light to the general conjecture of \cite{Tr}.

Recently,  Chanillo-Helffer-Laptev (\cite{CHL})
used Lidskii's Theorem to prove the existence of  nonlinear eigenvalues
of (\ref{p1})  permitting to recover some known results for  $n=1$
and to  give  new examples in dimension $n \ge 2$.

The proof of Theorem \ref{neven}  is based on the semiclassical analysis combined 
with the ideas of Chanillo-Helffer-Laptev (\cite{CHL}).
 We follow closely the reduction of Chanillo-Helffer-Laptev which is 
recalled in Section \ref{Section2}.
In this approach via Lidskii's theorem, the existence of a nonlinear 
eigenvalue of (\ref{p1}) is reduced to prove,  for some $k$, the non vanishing 
of the trace of the k-th power of a linear operator $\vD$:
\[
{\rm \;Tr \;} \vD^k \neq 0.
\]
Our approach is to introduce artificially  a semiclassical parameter $h$ 
in this operator (see Section \ref{Section3}).  Then, the existence of a nonlinear 
eigenvalue
of (\ref{p1}) is reduced to prove that for some $k$ and for some $h$,
\[
{\rm \;Tr \;} \vD(h)^k \neq 0,
\]
where $\vD(h)$ is an $h$-pseudo-differential operator. In other words,
 one can say that, while  Chanillo-Helffer-Laptev  try to
prove that, for some $k$, 
\begin{equation}
\sum_j \lambda_j^{-k}  \neq 0\;,
\end{equation}
where $\lambda_j$ are the  nonlinear eigenvalues of (\ref{p1}) (which can not be real under 
our assumptions), we 
want to find here  $h$ and $k$ such that
\begin{equation}
\sum_j (h^{m/(m+1)}\lambda_j+1)^{-k}  \neq 0.
\end{equation}
 It is standard (see \cite{Ro}) that ${\rm \;Tr \;} \vD(h)^k$ has a complete semiclassical expansion when $h \to 0$, 
which is theoretically computable. The question becomes now to find a 
non-zero term in these expansions which we call the semiclassical 
criterion.
The computation of the leading term given in Section \ref{Section4} gives the complete  
answer for $n$ even. 
 Although we do not have a general result for all $n$ odd, we believe 
that the semiclassical criterion  given below can be used to show that
 for each given value of $n$ odd,  (\ref{p1}) has at least one non 
trivial solution.
 The introduction of a semiclassical parameter $h$ permits us to overcome
the difficulty related to the non commutativity of operators encountered 
in (\cite{CHL}).

\section{Chanillo-Helffer-Laptev's approach}\label{Section2}

We rewrite $H(x, D_x, \lambda)$ in the form
\begin{equation}
H(x, D_x, \lambda) = L - 2 \lambda M + \lambda^2\;,
\end{equation}
with
\begin{equation}
L = - \Delta + P(x)^2\;,\; M = P(x)\;.
\end{equation}
 The operator $L$
 is invertible and its inverse is a pseudo-differential operator (see 
\cite{Ro}).
It is also easy to give sufficient conditions for determining whether
 the operator
\begin{equation}
A:= L^{-1}
\end{equation}
 belongs to a given Schatten class (see the appendix in \cite{CHL}).
Then the initial problem is reduced to the spectral analysis of
\begin{equation}\label{pbnl0}
(I - 2 \lambda B + \lambda^2 A) u = 0\;.
\end{equation}
with
\begin{equation}
B = A^{1/2} P A^{1/2}
\end{equation}

Chanillo-Helffer-Laptev (\cite{CHL}) used a rather standard  approach  
to transform
this nonlinear
spectral problem
to a linear ones  and then apply  Lidskii's Theorem to prove the 
existence
of a non trivial solution.

We first recall the reduction to the linear spectral problem.
It is easy to see that it is  enough to show that the operator $\mathcal 
D$ defined
 by
\begin{equation}\label{2.6}
\mathcal D := \left( \begin{array}{cc}
2 B & A^{1/2} \\
- A^{1/2}& 0
\end{array}
\right )
\end{equation}
has a non zero eigenvalue $\mu$. The first component
 of the corresponding eigenvector is an eigenvector of the problem (\ref{p1})
 with $\mu = \frac 1 \lambda$.

If $B$ and $A$ are compact, $\mathcal D$ is compact but the main
difficulty is
 that $\mathcal D$ is not selfadjoint.
The  Lidskii Theorem  says that

\begin{theorem}\label{thlidsk}.\\
Let $\mathcal C$ be a trace class operator. Let $\lambda_j (\mathcal
C)$ denote the sequence of non zero
eigenvalues of $\mathcal C$.  Then
$$\sum_j \lambda_j (\mathcal C) = {\rm \;Tr \;} \mathcal C\;.
$$
Here the eigenvalues are repeated according to their algebraic
multipicity.
\end{theorem}
By  Lidskii's Theorem, if ${\rm \; Tr \;} \mathcal C \neq 0$, ${\mathcal C}$ has at least
one non zero eigenvalue.

\begin{corollary} {\em (Chanillo-Helffer-Laptev \cite{CHL})}.\\
Let $\mathcal D$ be a compact operator. Assume that there exists $k\in
{\bf N}^*$
such that  $  \mathcal D ^k$  is of trace class and
$ {\rm \;Tr \;} \mathcal D ^k\neq 0.
$
Then (\ref{p1}) has at least one non trivial solution.
\end{corollary}

Chanillo-Helffer-Laptev used this criterion in the case $k=2,3,4$ and
obtained a family of interesting results in non analytic hypoellipticity.
The non commutativity between $A$ and $B$ makes it hard to
apply their method with higher rang criteria (although some results can 
be obtained in the same vein for $k=6,8$).

\section{The semiclassical criterion }\label{Section3}

As explained in the introduction,
our  goal is to introduce a semiclassical parameter $h$ in the operators 
$A$ and $B$
and to apply the theory of $h$-pseudo-differential operators
to give a complete semiclassical asymptotic expansion of the trace 
$\vD(h)^k$.
We are then led to find  conditions under which the leading term is non
zero.

Initially, the family of operators to study is of the form
\begin{equation}
H(x, D_x, \lambda) = - \Delta + (P(x) -\lambda )^2\;,
\end{equation}
where $x\mapsto P(x)$ is a real polynomial on $\mathbb  R^n$
 of order $m \ge 2$.  Write $P$ in the form
 $$
 P = P_m + P_{m-1} + \cdots +   P_0\;,
 $$
 where $P_j$ is a homogeneous polynomial of degree $j$. We assume that
 $P_m$ is elliptic on $\bR^n$:
 \begin{equation}\label{Pm}
 P_m(x) \neq 0, \quad \forall x \neq 0.
 \end{equation}
To be definite, we assume $P_m(x) > 0, \; \forall x \neq 0$.
The case $P_m(x) <  0, \;\forall x \neq 0$ is similar.
This condition on $P_m$ requires $m$ to be even and thus excludes
polynomials of odd degree when $n=1$.

To introduce the semiclassical parameter $h$, we make the dilation
$x \to \tau x$  and set $\lambda = (\lambda' -1) \tau^m$ (One takes
$\lambda =  (\lambda'+1) \tau^m$
if $P_m$ is negative).
Let
\begin{equation}
H(x, hD; \lambda', h) = - h^2\Delta + (  P(x,h) -\lambda' )^2\;,
\end{equation}
where $h =\frac{1}{\tau^{m+1}}$ and
\begin{equation}\label{expa}
P(x,h) = (P_m(x) +1) + h^{1/(m+ 1)}P_{m-1}(x)
+ h^{2/(m+ 1)}P_{m-2}(x)  + \cdots + h^{m/(m+ 1)}P_{0}\;.
\end{equation}
Then,
the initial problem
\begin{equation}\label{initial}
H(x, D; \lambda)v = 0
\end{equation}
has a non trivial solution $(\lambda, v)$ if and only if
\begin{equation}
H(x, hD; \lambda', h)u = 0
\end{equation}
has a non trivial solution $(\lambda', u)$ for some $h>0$.
Let us remark that
$$
-h^2 \Delta  + (P_m +1)^2 \ge 1.
$$ 
 One can then prove that for
$h>0$ small enough, $-h^2 \Delta + P(x,h)^2 > 1/2$ and therefore  is
invertible.\\

More generally, let us consider  the nonlinear eigenvalue problem 
\begin{equation} \label{h}
(-h^2 \Delta + (Q(x,h))^2 + ( P(x,h)  - z)^2)u =0,
\end{equation}
where $P(x,h)$ is of the above form and
\begin{equation}\label{expb}
Q(x,h) = Q_m(x)  + h^{1/(m+ 1)}Q_{m-1}(x)
+ h^{2/(m+ 1)}Q_{m-2}(x)  + \cdots + h^{m/(m+ 1)}Q_{0}\;,
\end{equation}
with $Q_j$ a homogenous polynomial of degree $j$. This kind of operators appears in the study of non analytic hypoellipticity for operators of the form (\ref{PQ}).
We assume that $P_m$ and $Q_m$ are real  and there exists $C>0$ such that
\begin{equation} \label{PQa}
 C^{-1} <x>^{2m} \le (P_m(x)+1)^2 + Q_m(x)^2  \le C <x>^{2m}, \quad x\in\bR^n.
\end{equation}
Here $<x> = ( 1 + |x|^2)^{1/2}$.
In the case $Q_m=0$, $P_m$ has to be an elliptic polynomial of degree $m$ so that
(\ref{PQa}) is satisfied.
We can show by the argument used above that $ -h^2\Delta + P(x,h)^2 + Q(x,h)^2$ is 
invertible for $h>0$ small.
We define the operators $A=A(h)$ and $B=B(h)$  in this new setting by :
\begin{equation}
A(h) = (-h^2\Delta + P(x,h)^2 + Q(x,h)^2)^{-1}, B(h) = A(h)^{1/2}
P(x,h)A(h)^{1/2}.
\end{equation}
For a temperate symbol, possibly $h$-dependent,  $K(x, \xi,h)$,
we denote by $K(x,hD, h)$ the $h$-pseudo-differential operators defined 
by
Weyl quantization
\begin{equation}
K(x,hD,h)u(x) =\frac{1}{(2\pi h)^n} \int \int e^{i(x-y)\cdot\xi/h} K(
\frac{x+y}{2}, \xi, h) u(y) \; dy d\xi\;,\; 
\end{equation}
$  \forall \;  u \in\mathcal S(\bR^n)$. By the  $h$-pseudo-differential
 calculus  \cite{Ro},
we deduce easily that $A(h)$ and $B(h)$ are $h$-pseudo-differential
operators
with symbols, satisfying  for any $N\in \nz^*$, 
\begin{eqnarray}
a(x,\xi;h) & = & \sum_{j = 0}^{(m+1)N} h^{\frac{j}{(m+1)}} a_j(x, \xi) +
h^{N + \frac{1}{(m+1)}}R_N(a,h) \\
b(x,\xi;h) & = & \sum_{j =0}^{(m+1)N} h^{\frac{j}{(m+1)}} b_j(x, \xi) + h^{N
+ \frac{1}{(m+1)}} R_N(b,h)\;.
\end{eqnarray}
Here the symbols 
\begin{equation}
\begin{array}{ll}
a_0 & =  (\xi^2 + (P_m +1)^2 + Q_m(x)^2)^{-1}\;,\\
b_0 &=  (\xi^2 + (P_m +1)^2+ Q_m(x)^2)^{-1}(P_m+1)\;,
\end{array}
\end{equation}
are, by definition, the $h$-principal symbols  of $A(h)$ and $B(h)$,
respectively.
If we denote by  $ S^{\rho}_{\phi,\varphi}$ the class of symbols defined as in
Robert \cite{Ro}, and if we introduce~:
$$
\begin{array}{ll}
 \rho_j &= ( 1 + \xi^2 + x^{2m})^{-1}
(1+x^2)^{-\frac{j}{2(m+1)}} (1+\xi^2)^{-\frac{j}{2(m+1)}},\\
\phi& = (1+x^2)^{1/2}, \\
\varphi &= (1+\xi^2)^{1/2}\;,
\end{array}
$$
then
$$
a_j \in S^{\rho_j}_{\phi,\varphi}, \quad b_j \in
S^{\rho_j<x>^m}_{\phi,\varphi}\;.
$$
Moreover the remainders  $R_N(a,h)$ and $R_N(b,h)$
are   respectively  a bounded family of symbols in $S^{\rho_{N(m+1)+1}}_{\phi,\varphi}$
and in $S^{\rho_{Nm+1}<x>^m}_{\phi,\varphi}$.

The Chanillo-Helffer-Laptev's $k$-criterion becomes now the following
\begin{lemma}\label{L3.1}.\\
  Let $\mathcal D (h)$ be
defined as in (\ref{2.6}) with $A$, $B$ replaced by $A(h)$, $B(h)$.
Let us assume there exists $k$ such that $\mathcal D(h)^k$ is an operator of 
trace class and that
\begin{equation}\label{kcritere}
 {\rm \;Tr \;} \mathcal D(h)^k \neq 0\;
\end{equation}
for some $h>0$. Then the nonlinear spectral problem (\ref{h}) has at
least one non trivial solution.
\end{lemma}

To apply this lemma,   we prove the following

\begin{theorem} \label{leth}.\\
Let us assume the condition  (\ref{PQa}) for $P$ and $Q$ with $m \ge 1$. Let $k > n (m+1)/m$, $n\ge 
1$.
Then $\mathcal D(h)^k$ is an operator  of trace class for
all $h>0$ sufficiently small.
We have, for any $N$,  the following asymptotic expansion
\begin{equation}\label{TrDk}
{\rm \;Tr \;} {\mathcal D}(h)^k = ( 2\pi h)^{-n} \{
\sum_{j=0}^{N(m+1)}
h^{\frac{j}{(m+1)}}H_{j; n,k}
+ \mathcal O(h^{N+ \frac{1}{(m+1)}}) \},
\end{equation}
when $ h \ar  0$. Here $H_{j;n,k}$ is independent of $h$ and $N$,  and can be 
computed from the symbol of $\vD(h)^k$. In particular,
\begin{equation}
H_{0;n, k} = \int_{\bR^{2n}} \tr \left( \sigma_k(x, \xi) \right) \; dx d\xi,
\end{equation}
$\sigma_k$ being the $h$-principal symbol of $\vD(h)^k$.
\end{theorem}
{\noindent \bf Proof}.\\
We note that $A(h)^{1/2}$ and $B(h)$ are  $h$-pseudo-differential operators
with symbol
in $S^{\rho_0^{1/2}}_{\phi, \varphi}$. Therefore $\mathcal D(h)^k$ is
an
 $h$-pseudo-differential operator  with matrix-valued symbol
in the class $S^{\rho_0^{\frac k2}}_{\phi, \varphi}$.
Since $ m \ge 1$, we have~:
$$
\rho_0^{\frac k2}\in L^1(\bR^{2n})\mbox{ if  } k > n(m+1)/m\;.
$$
Consequently,  $\mathcal D(h)^k$ is a trace class operator when
$k >n(m+1)/m$. If we denote by $\sigma_k(x, \xi;h)$ the total symbol
of $\vD(h)^k$, it  has a complete semiclassical 
expansion
similar to $a(h)$ beginning with $\sigma_k$, the $h$-principal symbol of
$\vD(h)^k$. The semiclassical expansion of the trace  follows from
the formula
$$
{\rm \;Tr \;} {\mathcal D}(h)^k = ( 2\pi h)^{-n} \int\int \tr \left(
\sigma_k(x, \xi; h) \right)\; d x d\xi.
$$
Here  $\tr$  denotes the trace of $2\times 2$ matrices. \qed

A consequence of Lemma \ref{L3.1}
 and  Theorem \ref{leth} is the following
\begin{corollary}{\em (The semiclassical criterion)}\label{cor3.3}.\\
Let (\ref{Pm}) be satisfied. If there exist $k\in\bN$ with $ k > (m+1)n/m$ and 
$j\in\bN$
such that $H_{j;n,k} \neq 0$,  then the nonlinear eigenvalue problem
(\ref{h}) has at least one non trivial solution $(z,u)$ for each $h>0$ sufficiently small.
 \end{corollary}

\section{An application of  the semiclassical criterion: the classical 
criterion}\label{Section4}

In this section, we apply the semiclassical criterion at the classical 
level, that is for 
$j=0$.

\begin{proposition}\label{Prop4.1}.\\
Assume that $Q_m =0$ and $P_m \ge 0$ is elliptic. Let $k >(m+1)n/m$.
Then one has
\begin{eqnarray}
H_{0;n,k} & =& 0, \quad \mbox{ if $n$ is odd, } \label{Hodd}\\
H_{0;n,k} & =& 2(-1)^\ell  C_n  \frac{(n-1)!}{(k-1)(k-2) \cdots(k-n)} 
C(P_m), \quad \mbox{ if $n= 2\ell$. }
\label{Heven}
\end{eqnarray}
Here $C_n$ is the volume of $S^{n-1}$ and
$$
 C(P_m) =  \int_{\bR^n_x} ( P_m(x) +1)^{-k} \; dx  >0 \;.
$$
\end{proposition}
In particular, we observe that~:
\begin{equation}
H_{0;n,k} \neq 0 \mbox{ when } k > (m+1)n/m 
\end{equation}
for all $n$ even.  As a consequence, we get Theorem \ref{neven}.\\

{\bf Proof of Proposition \ref{Prop4.1}}.\\
 The condition (\ref{PQa}) is satisfied. We compute 
 \begin{equation}
   \int\int \tr \left(\sigma_k(x, \xi)\right)\; dx d\xi.
 \end{equation}
where $\sigma_k$ is the $h$-principal symbol  of
 $\vD(h)^k$.
 Although the symbolic calculus for matrix-valued
$h$-pseudo-differential
 operators is complicated, the $h$-principal symbol of
 $\vD(h)^k$ can be easily computed.
 Since the $h$-principal symbol of $\mathcal D(h)$
 is
 \[ \left(\begin{array}{cc}
 2b_0 & a_0^{1/2} \\
 -a_0^{1/2} & 0
 \end{array}\right)
 \]
 the $h$-principal symbol of
 $\vD(h)^k$ is
 \[
 \sigma_k= \left(\begin{array}{cc}
 2b_0 & a_0^{1/2} \\
 -a_0^{1/2} & 0
 \end{array}\right)^k.
 \]
  Therefore,
 \begin{equation}
 \mbox{ tr } \sigma_k = ( b_0 + \sqrt{b_0^2-a_0} \; )^k + ( b_0 -
\sqrt{b_0^2-a_0}\;)^k\;.
 \end{equation}
 We recall that
\[ a_0=  (\xi^2 + (P_m +1)^2)^{-1}, \quad
b_0 =  (\xi^2 + (P_m +1)^2)^{-1}(P_m+1).
\]
By the change of variables $\xi \to (P_m +1)\eta$,
we obtain that
\begin{equation}
 \int\int \tr (\sigma_k)(x, \xi)\; dx d\xi = C(P_m)
\int_{\bR^n_\eta}( 1 +\eta^2)^{-k} ( (1 + i |\eta|)^k + (1 -i |\eta|)^k)
\; d\eta.
\end{equation}
This shows that
\begin{equation}
H_{0;n,k} = 2 C_n  C(P_m) \Re \int_0^\infty ( 1 + r^2)^{-k} (1+i r)^k r^{n-1} \;
dr, \quad k > (m+1)n /m,
\end{equation}
$C_n$ being the volume of $S^{n-1}$.\\

Let
\[
L(n,k) = \int_0^\infty ( 1 + r^2)^{-k} (1+i r)^k r^{n-1} \; dr 
 = \int_0^\infty  (1-i r)^{-k} r^{n-1} \; dr.
\]
An integration by parts gives~:
\begin{eqnarray}
L(n,k) &=& i\frac{n-1}{k-1}L(n-1, k-1)\\
&  = &i^{n-2}\frac{(n-1)(n-2)\cdots
2}{(k-1)(k-2) \cdots (k-n+2)}L(2, k-n+2) \;. \nonumber
\end{eqnarray}
Since
$$
L(2,j) = -\frac{1}{(j-1)(j-2)}\;,
$$
we have
\begin{equation}
H_{0;n,k} = 2 C_n\Re \{ i^{n}\frac{(n-1)!}{(k-1)(k-2) \cdots
(k-n+1)(k-n)}\} C(P_m).
\end{equation}
This proves (\ref{Hodd}) and (\ref{Heven}). \qed

Using  the classical criterion, we can  construct a family of examples when $n$ is odd. 

\begin{corollary}.\\
Let $n =n_1 + n_2$ with $n_1$ even and $n_2 \ge 1$. Let $P$, $R$  be real elliptic  
homogeneous polynomials\footnote{The important property is only
 that~: $- \Delta_{x''} + R(x'')^2$ on $\bR^{n_2}$ has an eigenvalue $\lambda_0$.}  of 
degree $m$  on $\bR^{n_1}$ and $\bR^{n_2}$, respectively.
For $x \in\bR^n$, set $x =(x', x'')\in\bR^{n_1}\times\bR^{n_2}$. Then,
the operator~:
$$
 -\Delta_x  +  \left(R(x'')D_{y_{1}}\right)^2 + 
 \left(P (x')  D_{y_{1}} - D_{y_{2}}\right)^2\;,
$$
is not  analytic hypoelliptic at $0$ in $\mathbb R^{n+2}$.
\end{corollary}
{\noindent \bf Proof}.\\
It suffices to show that the nonlinear eigenvalue problem
\[
\left( \sum_{j=1}^{n}  D_{x_j}^2  +  \left(R(x'')\right)^2 + 
 \left(P (x')   - z\right)^2 \right) u =0
\]
has a non trivial solution $(z, u)$. We can look for $v(x')$ satisfying, for some $z$,
\[
\left( - \Delta_{x'}  +  \lambda_0 +  \left(P (x')   - z'\right)^2 \right) v =0
\]
where $\lambda_0>0$ is an eigenvalue of $-\Delta_{x''} + R(x'')^2$. 
The corresponding semiclassical operator is  $-h^2 \Delta_{x'} + \lambda_0 h^{2m/(m+1)} + 
\left(P (x')   - z\right)^2 $ on $\bR^{n_1}$.
We can apply Proposition \ref{Prop4.1} for $Q = Q_0 = \lambda_0$  to conclude.
\qed

\begin{remark}.\\
One can play with $\lambda_0$. We can indeed consider a sequence $ \lambda_n$
 of eigenvalues of $-\Delta_{x''} + R(x'')^2$ tending to $+\infty$ and to
 choose a corresponding sequence $h_n$ of $h$'s tending to $0$,  such that 
$\lambda_n  h_n^{\frac{2m}{m+1}}=1$.
This could permit to relax the condition of ellipticity of $P$, particularly
 if we add a $Q$ like in (\ref{h}).
\end{remark}

Let us now consider the more general case $Q_m \neq 0$. For  $(P_m, Q_m)$ satisfying 
(\ref{PQa}), we can check that the classical criterion still works in some cases. Let us remark that
 \begin{equation}
 \mbox{ tr } \sigma_k = ( b_0 + \sqrt{b_0^2-a_0} \; )^k + ( b_0 -
\sqrt{b_0^2-a_0}\;)^k 
 \end{equation}
with now 
\[ a_0=  (\xi^2 + (P_m +1)^2 + Q_m^2 )^{-1}, \quad
b_0 =  (\xi^2 + (P_m +1)^2+ Q_m^2)^{-1}(P_m+1).
\]
Let 
\begin{equation}\label{tau}
T = ( ( P_m + 1)^2 + Q_m^2)^{1/2},\quad  \tau_1 = (P_m+1)/T, \quad \tau_2 = |Q_m|/T.
\end{equation}
By the change of variables $\xi \to T\eta$,
we obtain that
\begin{equation}\label{4.12}
 \int \tr (\sigma_k)(x, \xi)\; d\xi = 2 T^{n-k} \Re
\int_{\bR^n_\eta}( 1 +\eta^2)^{-k}  (\tau_1 + i (\tau_2^2 +\eta^2)^{1/2} )^k 
\; d\eta.
\end{equation}

\paragraph{The case $n=1$}.\\
This case  can be also thoroughly analyzed  by the method of \cite{CHL}.
When $n=1$,  we obtain  here for $k=3$,
\begin{equation}\label{H13}
H_{0;1,3} = -\frac{3 \pi}{2} \int_{\bR_x} \tau_1\tau_2^2 T^{-2} dx. 
\end{equation}
It is clear that if $P_m \ge 0$ (this requires $m$ to be even) and  $Q_m\neq 0$, 
then, $\tau_1 \ge 0$ and $H_{0;1,3} <0$.
If $P_m$ changes sign, $H_{0; 1, 3}$ may be vanishing. 
\\

\paragraph{The case $n \geq 2$}.\\
When $ n \ge 2$, $ k > (m+1)n/m$,
\begin{equation}
H_{0;n,k} = 2 C_n  \Re \int_{\bR^n_x} T(x)^{n-k} \int_0^\infty ( 1 + r^2)^{-k} (\tau_1 +i 
(\tau_2^2 + r^2)^{1/2})^k r^{n-1} \;dr.
\end{equation}
Set
\begin{equation}
C_{n,k} =  \int_0^\infty ( 1 + r^2)^{-k} (\tau_1 +i 
(\tau_2^2 + r^2)^{1/2})^k r^{n-1} \; dr. 
\end{equation}
By a change of variable, we obtain
\begin{equation*}
C_{n,k} = \int_{\tau_2}^\infty (\tau_1 -it)^{-k} (t^2- \tau_2^2)^{(n-2)/2} t \; dt.
\end{equation*}

\paragraph{The subcase $n=2$}.\\
In the case $n =2$,  it is easy to see that
\begin{equation}\label{4.16}
C_{2,k} =  \frac{(\tau_1 + i\tau_2)^{k-2}}{k-2} - \frac{\tau_1 (\tau_1 + 
i\tau_2)^{k-1}}{k-1}.
\end{equation}
Since $\tau_1^2 + \tau_2 ^2 =1$, we can see that
\begin{eqnarray}
\Re C_{2,3} & =& \tau_1 ( 1 + 2 \tau_2^2)/2  \ge 0,   \mbox{ if } \tau_1 \ge 0 \\
\Re C_{2,4} & =&-\frac 12 + 2\tau_1^2 - \frac 43 \tau_1^4.
\end{eqnarray}
This shows  $H_{0;2,3} >0$ for $n=2$ when $P_m \ge 0$. The classical criterion can be used 
in this case for $m>3$ when $n=2$.  \\

To study the case $n=2$ and $m =2$, we need information on  the sign of $H_{0;2,4}$ which 
depends on  the relation between $P_2$ and $Q_2$. An elementary computation 
 gives~: 
\begin{lemma}\label{casn=2}.\\
\begin{enumerate}
\item If $\tau_1 (x_1,x_2)^2 \le  (3 -\sqrt{3})/4$ for all $(x_1,x_2)$, then  $H_{0;2,4} 
>0$.
\item If $\tau_1 (x_1,x_2)^2 \ge  (3 -\sqrt{3})/4$ for all $(x_1,x_2)$,  then  $H_{0;2,4} 
<0$.
\end{enumerate}
\end{lemma}

\paragraph{The subcase $n=3$}.\\
When $n=3$, $k=5$, we have 
\begin{equation}
H_{0;3,5} = \frac{5 \pi}{8} \int_{\bR_x^3} \tau_2^4 \tau_1 T^{-2} dx >0, 
\end{equation}
so long as $Q_m$ is not identically zero and $P_m \ge 0$.\\

\paragraph{Some applications}.\\
We obtain the following~:
\begin{proposition} \label{Prop4.3}.\\
(a). Let $n=1$. The operator
\[
P_m(x, D) := D_x^2 + (D_s - c x^m D_t)^2 + (x^m D_t)^2
\] 
 is not analytic hypoelliptic  at $0\in \bR^3$ for all  $m \ge 2$ even and 
$c\in\bR$. 

(b). Let $n=2$. Let $P$, $Q$ be real homogenous polynomials of degree $m \ge 2$
satisfying (\ref{PQa}) and $P \ge 0$.  When $m =2$, we assume in addition that one of the 
above conditions on $\tau_1$ is satisfied. Then the operator
\begin{equation}
 D_x^2 + D_y^2  + (D_s- P(x,y)D_t)^2  + (Q(x,y)D_t)^2, \quad (x,y)\in\bR^3,
\end{equation}
is not analytic hypoelliptic at $0$ in $\bR^4$.

(c). Let $n=3$. Let $P$, $Q$ be real homogenous polynomials on $\bR^3$ of degree $m \ge 2$ 
satisfying (\ref{PQa}). Assume that $Q \neq 0, P \ge 0$. Then the operator
\begin{equation}
 - \Delta_x  + (D_s- P(x)D_t)^2  + (Q(x)D_t)^2
\end{equation}
is not analytic hypoelliptic at $0$ in $\bR^5$.
\end{proposition}

We believe that the condition on $\tau_1$ in (b) is technical.
As examples of $(P, Q)$ satisfying the conditions of (b) of Proposition 
\ref{Prop4.3},
we can take  $P(x,y) = (x^2 + y^2)^{\ell}$,  $Q(x,y) = (xy)^{\ell}$, 
$\forall \; \ell  \ge 1$, because in the case  $\ell=1$, one can check that $\tau_1(x,y)^2 
\ge 
4/5 > (3-\sqrt{3})/4$. We shall come back at the analysis of this example
corresponding to $\ell=1$  in the next section.

\begin{remark}.\\
Let us briefly compare the results  obtained here and those of 
\cite{CHL}.
Apparently, in the case $Q =0$, the ``classical'' criterion does not produce any result for 
$n \ge 1$ odd. But the semiclassical criteria  may still work by 
considering a  coefficient of  higher
order of the expansion of the trace. 
As indicated by J.~Sj\"ostrand, a similar approach is
 used by L.~Nedelec in \cite{Ne} for getting a lower bound for the number of resonances of an 
$h$-pseudodifferential system. The same condition on the dimension appears.
 The ``quantum'' criterion given in \cite{CHL} works for
 $n=1,2,3$ but with a stronger condition on $m$ when $n>1$. Moreover
it has been  observed in Remark 4.4  in \cite{CHL} that  the condition of ellipticity
 of $P$ can be replaced by a weaker condition.
The last point is that the homogeneity of $P$ plays an  important role
 in the dilation argument of \cite{CHL}, while in the semiclassical 
approach,  the lower order parts  can be included. This appears to be  useful in the 
dimension reduction.
\end{remark}

\section{Comparison with M\'etivier's results}\label{Section5}
In this section, we would like to analyze the links between our results
 and the previous results obtained by G.~M\'etivier \cite{Met3, Met, Met5}.
\subsection{A first family of examples}
Let us first consider the operator on $\bR^{n+2}$~:
$$
H(X, D_X)=- \Delta + (P(x) D_{x_{n+1}} - D_{x_{n+2}})^2\;,
$$
where $X=(x,x_{n+1},x_{n+2})$,  $P$ is an homogeneous positive elliptic polynomial of degree 
$m \geq 2$ on $\bR^n$.\\

If we take the ``microlocal spirit'' we observe that $H$ is an operator with
 double characteristics, whose principal symbol is the function~:
$$
(T^* \bR^{n+2}\setminus 0 )  \ni (X, \Xi) \mapsto |\xi|^2 +
 (P(x) \xi_{n+1} - \xi_{n+2})^2\;.
$$
The symbol vanishes exactly at order $2$ on the submanifold
$$
\Sigma = \{ (X,\Xi)\;|\; \xi = 0\;,\;P(x) \xi_{n+1} - \xi_{n+2} =0\;,\; \xi_{n+1}\neq 0\}\;.
$$
This submanifold is of codimension $n+1$. Let us now analyze the ``symplecticity''
 of $\Sigma$. We recall that $\Sigma$ is said to be symplectic if the
 restriction  to $\Sigma$ of the canonical $2$-form is non degenerate.
An easy way for verifying the symplecticity is to consider the 
$(n+1)\times (n+1)$ matrix
${u_i,u_j}$ where $u_i (X, \Xi)= \xi_i$ for $i=1,\cdots,n$ and $u_{n+1}(X,\Xi)=
P(x) \xi_{n+1} - \xi_{n+2}$ and to show that it is not degenerate.
An immediate computation shows that its rank at a given point is $2$ if $\nabla P\neq 0$
 and $0$ if $\nabla P =0$.\\
 
When $P$ is elliptic and homogeneous, we get that the rank is constant outside $0$ and equal 
to $2$.
There are two cases~:
\begin{enumerate}
\item When $n=1$, we get that $\Sigma$ is symplectic except at the points
 of $\Sigma$ such that $x=0$. The result of Tr\`eves-Tartakoff-M\'etivier-Sj\"ostrand 
\cite{Tar, Tr1, Met4, Sj}  shows
 that the operator is microlocally analytic hypoelliptic  outside of $\Sigma$ (ellipticity) 
and in the neighborhood of the points of $\Sigma$ such that $x\neq 0$.
In this case, the operator is not analytic hypoelliptic at any point 
$(0,x_{n+1}, x_{n+2})$.
\item When $n >1$, M\'etivier's result\footnote{Note that the operator is hypoelliptic with 
loss of one derivative in $\bR^{n+2} \setminus \{x=0\}$.}  \cite{Met3, Met} gives that the 
operator $H$
 is not analytic hypoelliptic in any open set in $\bR^{n+2}$. What we show here is the more 
precise 
result that $P$ is not  analytic hypoelliptic at any point $(0,x_{n+1}, x_{n+2})$ which is a 
finer property. See the introduction of \cite{Met}
 for the discussion between the definitions of analytic  hypoellipticity and 
germ-hypoanalyticity (analytic hypoellipticity in a neighborhood of a point
 and analytic hypoellipticity at a point).
\end{enumerate}

\subsection{A new class of non analytic hypoelliptic operators}
We now show that maybe more interesting examples can be treated 
 when considering  the more general class~:
$$
H(X, D_X)=- \Delta + (P(x) D_{x_{n+1}} - D_{x_{n+2}})^2 + Q(x)^2 D_{x_{n+1}}^2\;,
$$
where $P$ and $Q$ are homogeneous polynomials of degree $m >1$ with $P \geq 0$
 and $P^2 + Q^2$ elliptic.
 When restricting $H$ to $x_{n+1}$ independent distributions, we get
 an analytic hypoelliptic operator on $\bR^{n+1}$:
 $$- \Delta_x + (P(x)^2 + Q(x)^2) D_{x_{n+2}}^2\;,
$$
by applying a theorem of Grushin (\cite{Gr}).
\\

We have seen in Section \ref{Section4}, that the ``classical'' criterion
 can give a result under some additional condition.
We will concentrate our analysis to the specific case when~:
$$
n=2\;,\; P(x) = x_1^2 + x_2^2\;,\; Q(x) = \alpha x_1 x_2\;,
$$
with $\alpha >0$.\\

Let us do the same microlocal analysis as in the previous subsection. The
 characteristic set $\Sigma$ is now defined as a union of two regular submanifolds of 
dimension $4$
 in $\bR^8 \setminus 0$~:
$$
\Sigma = \Sigma_1 \cup \Sigma_2 \;,
$$
with
$$
\Sigma_j = \{ \xi_1=0 \;,\; \xi_2 =0 \;,\; \xi_4 = (x_1^2 + x_2^2) \xi_3\;,\; x_j =0, 
\xi_3\neq 0\}\;.
$$
Moreover $\Sigma_j$ is symplectic outside $\Sigma_1 \cap \Sigma_2$
 and not symplectic at $\Sigma_1\cap \Sigma_2$.\\
 
Outside $\Sigma_1\cap \Sigma_2$, observing that the symbol
 of $H$ vanishes exactly at order $2$ on $\Sigma$ there,  we get again from \cite{Tar, Tr1, 
Met4, Sj}
 that $H$ is microlocally analytic hypoelliptic.
M\'etivier's criterion of non-analytic-hypoellipticity can not be applied at the
 points $(0,0, x_3, x_4)$ (the operator is indeed not hypoelliptic with loss of one derivative) and 
it is interesting to see what is obtained 
 through our approach.
\begin{proposition}\label{prop5.1}.\\
For $\alpha >0$ small enough, the operator $H(X,D_X)$ is not analytic 
hypoelliptic at
any point  $(0,x_3,x_4)$.
\end{proposition}
{\noindent \bf Proof}.\\
We just apply the criterion of the previous section (Proposition \ref{Prop4.3}, second case) 
and the discussion following the statement.

When $\alpha$ is small enough (at least $0<\alpha \leq 1$), we observe that 
the quantity 
 $ (P_m+1)^2 / \left((P_m+1)^2 + Q_m^2\right)$ is sufficiently 
 near $1$. In particular, 
the second condition on $\tau_1$ in Lemma \ref{casn=2} is satisfied.

\begin{remark}.\\
As explained to us by M.~Christ, it is possible to prove, 
 for rather large classes of models depending analytically on  an additional parameter $\alpha$, that  some associated Fredholm determinant is analytic in $\alpha$.  One can  conclude that, if the operator is not analytic hypoelliptic
for some value of $\alpha$, then it is not  analytic
hypoelliptic for generic values of $\alpha$. We refer to  \cite{Chr3}, Proposition 5.2
 for the argument.
In the particular case of the above Proposition \ref{prop5.1}, we can present the argument in the following way.  Let $H_{0,2,4}(\alpha)$ denote $H_{0,2,4}$ defined as in Theorem \ref{leth} with $P = x_1^2 + x_2^2$ and  $Q = \alpha x_1 x_2$. Using (\ref{4.12}) and (\ref{4.16}), one can check that $H_{0,2,4}(\alpha)$ is real analytic in $\alpha>0$. The proof of Proposition  \ref{prop5.1}  shows that
$H_{0,2,4}(\alpha) \neq 0$ for $\alpha >0$ small. Thus, $H_{0;2,4}(\alpha) \neq 0$ for all $\alpha>0$, except a discrete set in $\bR_+$ and therefore  Proposition  \ref{prop5.1}  remains true in this case. This argument can be used  for  more general cases. It is indeed  easy to check in many 
situations  the analyticity of  $H_{0,n,k}$ with respect to the parameter. 
\end{remark}


\begin{thebibliography}{99}


\bibitem{Cha} S. Chanillo~:
\newblock Kirillov theory, Tr\`eves strata, Schr\"odinger equations and
analytic
hypoellipticity of sums of squares.
\newblock Preprint August 2001,  http://arxiv.org/pdf/math.AP/0107106).


\bibitem{CHL} S. Chanillo, B. Helffer and A. Laptev~:
\newblock Non linear eigenvalues and analytic hypoellipticity. 
\newblock Preprint Institut Mittag-Leffler and 
http://arxiv.org/pdf/math.AP/0211308.
\newblock To appear in Journal of Functional Analysis 2003.


\bibitem{Chr1} M. Christ~:
\newblock Some non-analytic-hypoelliptic sums of squares of vector
fields.
\newblock Bull. A.M.S 16, p.~137-140  (1992).

\bibitem{Chr2} M. Christ~:
\newblock Analytic hypoellipticity, representations of nilpotent
groups, and a non-linear eigenvalue problem.
\newblock Duke Math. J. 72, p.~595-639  (1993).


\bibitem{Chr3} M. Christ~:
\newblock Analytic hypoellipticity in dimension two. 
\newblock  MSRI Preprint, No. 1996-009 (1996).


\bibitem{Chr4} M. Christ~:
\newblock The Szeg\"o projection need not preserve global analyticity.
\newblock Ann. of Math., 301-330, 143(1996).


\bibitem{Gr} V.V.~Grushin~:
\newblock On a class of elliptic pseudodifferential 
operators degenerate on a submanifold.
\newblock Math. USSR Sbornik, Vol. 13 (1971), No 2, p.~155-185.

\bibitem{HH} N. Hanges, A.A.  Himonas~:
\newblock Non-analytic hypoellipticity in the presence of
symplecticity.
\newblock
Proc. Am. Math. Soc. 126, n$^\circ$2, p.~405-409 (1998). 

\bibitem{He} B. Helffer~:
\newblock  Conditions n\'ecessaires d'hypoanalyticit\'e pour des
op\'erateurs invariants \`a gauche sur un groupe nilpotent gradu\'e.
\newblock Journal of differential Equations, Vol.~44, n$^\circ$3, 
p.~460-481
(1982).

\bibitem{Ma} A.S. Markus~:
\newblock Introduction to the spectral theory of polynomial operator 
pencils.
\newblock Vol. 71, Translations of mathematical monographs. American 
Mathematical Society.

\bibitem{Met3} G. M\'etivier~:
\newblock Une classe d'op\'erateurs non hypoelliptiques analytiques.
\newblock S\'eminaire Goulaouic-Schwartz 1978-1979, Expos\'e XII.

\bibitem{Met4} G. M\'etivier~:
\newblock Analytic hypoellipticity for operators with multiple
 characteristics.
\newblock Comm. in PDE 6, p.~1-90 (1980).

\bibitem{Met} G. M\'etivier~:
\newblock Une classe d'op\'erateurs non-hypoelliptiques analytiques.
\newblock Indiana Univ. Math. J., Vol.~29, p.~823-860 (1980).

\bibitem{Met5} G. M\'etivier~: Non hypoellipticit\'e analytique pour des
op\'erateurs \`a caract\'eristiques doubles, 
\newblock S\'eminaire Goulaouic-Meyer-Schwartz 1981-1982, Expos\'e XII.

\bibitem {Ne} L.~Nedelec~:
\newblock Existence of resonances for matrix Schr\"odinger operators.
\newblock To  appear in Asymptotic Analysis (2003).

\bibitem{PhRo}  Pham The Lai, D. Robert~:
\newblock Sur un probl\`eme aux valeurs propres non lin\'eaire,
\newblock Israel Journal of Math. 36, p.~169-186 (1980).

\bibitem{Ro} D.~Robert~:
\newblock {\it Autour de l'approximation semi-classique}.
\newblock Progress in Mathematics n$^0$ 68, Birkh\"auser (1987).

\bibitem{Sim} B. Simon~:
\newblock Trace ideals and their applications.
\newblock London Mathematical Society. Lecture Note Series 35. Cambridge 
University
 Press (1979).


\bibitem{Sj} J. Sj\"ostrand~:
\newblock Analytic wavefront sets and operators with multiple
characteristics.
\newblock Hokkaido Mathematical Journal, Vol.~12, p.~393-433 (1983).


\bibitem{Tar} D. Tartakoff~:
\newblock On the local real analyticity of solutions to $\Box _b$ and the 
${\bar \pa}$-Neumann problem.
\newblock Acta Math. 145, p.~117-204 (1980).


\bibitem{Tr1} F. Tr\`eves~:
\newblock Analytic hypoellipticity of a class of pseudo-differential
operators
 with double characteristics and applications to the ${\bar
\pa}$-Neumann problem.
\newblock Comm. in PDE 3, p.~476-642  (1978).


\bibitem{Tr} F. Tr\`eves~:
\newblock Symplectic geometry and analytic hypoellipticity.
\newblock Differential equations: La Pietra 1996 (Florence), p.~201-219,
\newblock Proc.
   Sympos. Pure Math., 65, Amer. Math. Soc., Providence, RI, 1999.

\end{thebibliography}
\end{document}